\documentclass{amsart}

\usepackage{amsmath,amsthm,amsfonts,amscd,amssymb}

\newtheorem{thm}{Theorem}[section]

\newtheorem{lem}[thm]{Lemma}

\theoremstyle{definition}

\theoremstyle{remark}
\newtheorem{rem}{Remark}[section]

\begin{document}

\title[Orthogonal Siegel's lemma for real polynomial spaces]{Integral orthogonal bases of small height for real polynomial spaces}
\author{Lenny Fukshansky}

\address{Department of Mathematics, 850 Columbia Avenue, Claremont McKenna College, Claremont, CA 91711}
\email{lenny@cmc.edu}

\begin{abstract} Let $\mathcal P_N(\mathbb R)$ be the space of all real polynomials in $N$ variables with the usual inner product $\left<\ ,\ \right>$ on it, given by integrating over the unit sphere. We start by deriving an explicit combinatorial formula for the bilinear form representing this inner product on the space of coefficient vectors of all polynomials in $\mathcal P_N(\mathbb R)$ of degree $\leq M$. We exhibit two applications of this formula. First, given a finite dimensional subspace $V$ of $\mathcal P_N(\mathbb R)$ defined over $\mathbb Q$, we  prove the existence of an orthogonal basis for $\left( V, \left<\ ,\ \right> \right)$, consisting of polynomials of small height with integer coefficients, providing an explicit bound on the height; this can be viewed as a version of Siegel's lemma for real polynomial inner product spaces. Secondly, we derive a criterion for a finite set of points on the unit sphere in $\mathbb R^N$ to be a spherical $M$-design.
\end{abstract}

\maketitle

\def\A{{\mathcal A}}
\def\AA{{\mathfrak A}}
\def\B{{\mathcal B}}
\def\C{{\mathcal C}}
\def\D{{\mathcal D}}
\def\F{{\mathcal F}}
\def\x{{\mathcal H}}
\def\I{{\mathcal I}}
\def\J{{\mathcal J}}
\def\K{{\mathcal K}}
\def\kk{{\mathfrak K}}
\def\L{{\mathcal L}}
\def\LL{{\mathfrak L}}
\def\M{{\mathcal M}}
\def\mm{{\mathfrak m}}
\def\MM{{\mathfrak M}}
\def\OO{{\mathfrak O}}
\def\R{{\mathcal R}}
\def\PNR{{\mathcal P_N(\real)}}
\def\PMNR{{\mathcal P^M_N(\real)}}
\def\PdNR{{\mathcal P^d_N(\real)}}
\def\s{{\mathcal S}}
\def\V{{\mathcal V}}
\def\X{{\mathcal X}}
\def\Y{{\mathcal Y}}
\def\Z{{\mathcal Z}}
\def\H{{\mathcal H}}
\def\cee{{\mathbb C}}
\def\pee{{\mathbb P}}
\def\que{{\mathbb Q}}
\def\real{{\mathbb R}}
\def\zed{{\mathbb Z}}
\def\aaa{{\mathbb A}}
\def\ff{{\mathbb F}}
\def\kk{{\mathfrak K}}
\def\qbar{{\overline{\mathbb Q}}}
\def\kbar{{\overline{K}}}
\def\ybar{{\overline{Y}}}
\def\kkbar{{\overline{\mathfrak K}}}
\def\ubar{{\overline{U}}}
\def\eps{{\varepsilon}}
\def\ahat{{\hat \alpha}}
\def\bhat{{\hat \beta}}
\def\gt{{\tilde \gamma}}
\def\h{{\tfrac12}}
\def\be{{\boldsymbol e}}
\def\bei{{\boldsymbol e_i}}
\def\bc{{\boldsymbol c}}
\def\bm{{\boldsymbol m}}
\def\bk{{\boldsymbol k}}
\def\bi{{\boldsymbol i}}
\def\bl{{\boldsymbol l}}
\def\bq{{\boldsymbol q}}
\def\bu{{\boldsymbol u}}
\def\bt{{\boldsymbol t}}
\def\bs{{\boldsymbol s}}
\def\bv{{\boldsymbol v}}
\def\bw{{\boldsymbol w}}
\def\bx{{\boldsymbol x}}
\def\bX{{\boldsymbol X}}
\def\bz{{\boldsymbol z}}
\def\bwy{{\boldsymbol y}}
\def\bY{{\boldsymbol Y}}
\def\bL{{\boldsymbol L}}
\def\ba{{\boldsymbol a}}
\def\bb{{\boldsymbol b}}
\def\bet{{\boldsymbol\eta}}
\def\bxi{{\boldsymbol\xi}}
\def\bo{{\boldsymbol 0}}
\def\bone{{\boldsymbol 1}}
\def\bol{{\boldsymbol 1}_L}
\def\ep{\varepsilon}
\def\p{\boldsymbol\varphi}
\def\q{\boldsymbol\psi}
\def\rank{\operatorname{rank}}
\def\aut{\operatorname{Aut}}
\def\lcm{\operatorname{lcm}}
\def\sgn{\operatorname{sgn}}
\def\spn{\operatorname{span}}
\def\md{\operatorname{mod}}
\def\Norm{\operatorname{Norm}}
\def\dim{\operatorname{dim}}
\def\det{\operatorname{det}}
\def\Vol{\operatorname{Vol}}
\def\rk{\operatorname{rk}}
\def\ord{\operatorname{ord}}
\def\ker{\operatorname{ker}}
\def\div{\operatorname{div}}
\def\Gal{\operatorname{Gal}}

\section{Introduction and notation}
\label{intro}

Siegel's lemma originated as an important combinatorial principle that a system of homogeneous linear equations with integer coefficients should have a nontrivial integral solution vector whose entries are comparable in size to the coefficients of the system. Although this observation was already made by Thue \cite{thue} in 1909, the first formal proof of such a result by an application of the pigeonhole principle appeared in the paper \cite{siegel} of Siegel in 1929. In its modern formulation, Siegel's lemma is a statement about the existence of a "short" basis for a vector space over a global field, where the size of the vectors is measured with respect to a {\it height function}, a standard tool of Diophantine geometry which generalizes the naive sup-norm over integers. A general result like this was first proved over number fields by Bombieri and Vaaler \cite{vaaler:siegel} in 1983 with further extensions by a variety of authors following in the consequent years. In particular, given a symmetric bilinear space one may ask for a short {\it orthogonal} basis over a fixed field or ring. Results of this nature were recently obtained in \cite{me:witt} and \cite{me:qbar_quad}, where the height of basis vectors in question were bounded in terms of the heights of the vector space and the coefficient vector of the bilinear form. One goal of the present note is to produce a similar result for real polynomial spaces. We start by setting up some notation.

Let $N \geq 2$ be an integer, and let us consider the algebra of all polynomials in $N$ variables with real coefficients
$$\PNR := \real[X_1,...,X_N]$$
as an infinite-dimensional real vector space with the standard inner product on it (see, for instance Chapter IV of \cite{stein:weiss}), given by
\begin{equation}
\label{inner:product}
\left< F,G \right> = \frac{1}{\alpha_N} \int_{\Sigma_{N-1}} F(\bx) G(\bx) d\bx,
\end{equation}
for each $F, G \in \PNR$, where $\Sigma_{N-1}$ is the unit sphere in $\real^N$, we integrate with respect to the usual Lebesgue measure on $\real^N$, and $\alpha_N$ is the generalized surface area of $\Sigma_{N-1}$:
\begin{equation}
\label{surface}
\alpha_N = \left\{ \begin{array}{ll}
\frac{(2 \pi)^{N/2}}{2 \times 4 \times \dots \times (N-2)} & \mbox{if $N$ is even} \\
\frac{2(2 \pi)^{(N-1)/2}}{1 \times 3 \times \dots \times (N-2)} & \mbox{if $N$ is odd}.
\end{array}
\right.
\end{equation}
Let $M \geq 1$ be an integer, and write
\begin{equation}
\label{MNM}
\M(M,N) = \left\{ \bm \in \zed_{\geq 0}^N : w(\bm) := \sum_{i=1}^N m_i \leq M \right\}
\end{equation}
for a set of multi-indexes, arranged in lexicographic order; we call $w(\bm)$ the {\it weight} of the index vector $\bm$. Define $\PMNR$ to be the space of all polynomials in $\PNR$ of degree $\leq M$, then each polynomial $F(X_1,...,X_N) \in \PMNR$ can be written as
\begin{equation}
\label{F_c}
F(\bX) = \sum_{\bm \in \M(M,N)} c_F(\bm) \bX^{\bm},
\end{equation}
where $\bX^{\bm} = X_1^{m_1}...X_N^{m_N}$, $c_F(\bm) \in \real$ for each $\bm \in \M(M,N)$; here and for the rest of the paper we adopt the convention that $X_i^0=1$, even when $X_i=0$. We also write $\bc_F = (c_F(\bm))_{\bm \in \M(M,N)} \in \real^{L(M,N)}$ for the vector of coefficients of $F$, where
\begin{equation}
\label{MNM.1}
L(M,N) := |\M(M,N)| = \dim_{\real} \PMNR = \sum_{k=0}^M \binom{N+k-1}{k}.
\end{equation}
An important tool we need is an explicit combinatorial formula for $\left< F, G \right>$ in terms of the coefficients of polynomials $F(\bX)$ and $G(\bX)$.
\smallskip

Let us fix $M \geq 1$, and let $L = L(M,N)$ as given by (\ref{MNM.1}). Let us also define the double factorial $m!!$ for any integer $m \geq -1$ to be
\begin{equation}
\label{factorial}
m!! = \left\{ \begin{array}{ll}
m(m-2)(m-4) \dots 5 \times 3 \times 1 & \mbox{if $m > 1$ is odd} \\
m(m-2)(m-4) \dots 6 \times 4 \times 2 & \mbox{if $m > 1$ is even} \\
1 & \mbox{if $m=-1,0,1$}.
\end{array}
\right.
\end{equation}
For each $\bm = (m_1,...,m_N) \in \M(M,N)$, define
\begin{equation}
\label{pm:def}
P(\bm) = \frac{\prod_{i=1}^N (2m_i-1)!! }{\prod_{k=1}^{w(\bm)} \left( N - 2 + 2k \right)},
\end{equation}
where an empty product $\prod_{k=1}^0$ is interpreted as equal to 1. For each $\ba \in \real^L$, we write $\ba = (a(\bm))_{\bm \in \M(M,N)}$. Let 
$$2\M(M,N)= \left\{ 2\bm : \bm \in \M(M,N) \right\},$$
and let 
$$E(M,N) = \left\{ (\bm_1,\bm_2) \in \M(M,N) \times \M(M,N) : \bm_1+\bm_2 \in 2\M(M,N) \right\}.$$
Define a bilinear form $\L_{M,N}\ : \real^L \times \real^L \longrightarrow \real$ by
\begin{equation}
\label{bilin:def}
\L_{M,N}(\ba,\bb) = \sum_{(\bm_1,\bm_2) \in E(M,N)} P\left( \frac{\bm_1+\bm_2}{2} \right) a(\bm_1)b(\bm_2),
\end{equation}
for each $(\ba,\bb) \in \real^L \times \real^L$, and let $\L_{M,N}(\ba):= \L_{M,N}(\ba,\ba)$ be the corresponding quadratic form. Notice in particular that if $(\bm_1,\bm_2) \in \M(M,N) \setminus E(M,N)$, then the coefficient of $\L_{M,N}(\ba,\bb)$ corresponding to the monomial $a(\bm_1) b(\bm_2)$ is equal to zero. Then we have the following result.
\smallskip

\begin{thm} \label{int_lemma} For each $F(\bX),G(\bX) \in \PMNR$,
\begin{equation}
\label{integral:evaluation}
\left< F,G \right> = \L_{M,N}(\bc_F, \bc_G),
\end{equation}
where $\bc_F,\bc_G$ are coefficient vectors of $F$ and $G$, respectively.
\end{thm}

\begin{rem} \label{norms} An immediate implication of  (\ref{integral:evaluation}) is that the quadratic form $\L_{M,N}(\ba)$ must be positive definite, since it is a norm form.
\end{rem}

We prove Theorem \ref{int_lemma} in section~\ref{integral}. Next, let $V$ be a $n$-dimensional subspace of $\PNR$ defined over $\que$, $n \geq 1$, and define {\it degree} of $V$ to be
$$d = \deg(V): = \min \left\{ M \in \zed : V \subseteq \PMNR \right\}.$$
Then $\left( V, \left<\ ,\ \right> \right)$ is a finite-dimensional $\que$-inner product space, so there must exist an orthogonal basis for $V$ consisting of polynomials with integer coefficients. We will prove the existence of a basis like this with each polynomial having relatively small height, which is a version of Siegel's lemma for polynomial spaces with additional orthogonality conditions. 

For each polynomial $F(\bx) \in \PNR$, define the height of $F$ by
\begin{equation}
\label{height_poly}
H(F) = |\bc_F| := \max_{\bm \in \M(M,N)} |c_F(\bm)|,
\end{equation}
where $M = \deg(F)$. Let $L = |\M(d,N)|$, and define an injective linear map $\varphi : V \to \real^L$ by $\varphi(F) = \bc_F$. Let $f_1,\dots,f_k$ be a basis for $V$, and let $C = (\bc_{f_1} \dots \bc_{f_k})$ be the corresponding $L \times k$ basis matrix for $\varphi(V)$. Then define the height of $V$ by
$$H(V) = D^{-1} \sqrt{|\det(C^t C)|},$$
where $D$ is the greatest common divisor of the determinants of all $k \times k$ minors of~$C$; $H(V)$ is well-defined, i.e. this definition does not depend on the choice of a basis for $V$. With this notation, we can now state our main result.

\begin{thm} \label{int_basis} Let $V$ be as above. Then there exists an orthogonal basis $g_1,\dots,g_n$ for $\left( V, \left<\ ,\ \right> \right)$ consisting of polynomials with integer coefficients so that
\begin{equation}
\label{main_bound}
\prod_{i=1}^n H(g_i) \leq L^{\frac{3n(n+1)}{2}} H(V)^n.
\end{equation}
\end{thm}

We prove Theorem \ref{int_basis} in section~\ref{basis}. Our main tools are Theorem \ref{int_lemma} and a version of Siegel's lemma due to Bombieri and Vaaler \cite{vaaler:siegel}. In fact, we prove a more general statement than Theorem \ref{int_basis}, Theorem \ref{gen}, where $\left<\ ,\ \right>$ is replaced by an arbitrary bilinear form defined on $V$. Finally, in section \ref{sphere} we derive an application of Theorem \ref{int_lemma} to spherical designs, presenting a necessary and sufficient criterion for a finite set of points on the unit sphere in $\real^N$ to be a spherical $M$-design (see Theorem \ref{design}).
\bigskip

\section{Proof of Theorem \ref{int_lemma}}
\label{integral}

Let us write $\M$ for $\M(M,N)$, $E$ for $E(M,N)$, and $\L$ for $\L_{M,N}$. First notice that
\begin{equation}
\label{in1}
\left< F,G \right> = \sum_{\bm_1 \in \M} \sum_{\bm_2 \in \M} c_F(\bm_1) c_G(\bm_2) S(\bm_1,\bm_2),
\end{equation}
where
\begin{equation}
\label{in2}
S(\bm_1,\bm_2) = \frac{1}{\alpha_N} \int_{\Sigma_{N-1}} \bx^{\bm_1+\bm_2} d\bx = \frac{1}{\alpha_N} \int_{\Sigma_{N-1}} \prod_{i=1}^N x_i^{\eps_i} d\bx,
\end{equation}
where the weight $w(\bm_1+\bm_2) = \sum_{i=1}^N \eps_i \leq 2M$, $\eps_i \in \zed_{\geq 0}$ for all $1 \leq i \leq N$. Consider a change to spherical coordinates (see, for instance, page 181 of \cite{fleming}) $0 \leq \theta_i \leq \pi$ for all $1 \leq i \leq N-2$, $0 \leq \theta_{N-1} \leq 2\pi$, given by 
\begin{equation}
\label{in3}
x_i = \cos\ \theta_i \prod_{j=1}^{i-1} \sin\ \theta_j,
\end{equation}
for all $1 \leq i \leq N-1$, and $x_N = \prod_{j=1}^{N-1} \sin\ \theta_j$. The Jacobian of this coordinate change is
\begin{equation}
\label{in4}
J=\prod_{i=1}^{N-2} \sin^{N-1-i} \theta_i.
\end{equation}

\noindent
Then
\begin{eqnarray}
\label{in5}
S(\bm_1,\bm_2) & = & \frac{1}{\alpha_N} \left( \prod_{i=1}^{N-2} \int_0^{\pi} \cos^{\eps_i} \theta_i\ \sin^{\beta_i} \theta_i\ d\theta_i \right) \times \nonumber \\
& \times & \int_0^{2\pi} \cos^{\eps_{N-1}} \theta_{N-1}\ \sin^{\beta_{N-1}} \theta_{N-1}\ d\theta_{N-1},
\end{eqnarray}
where $\beta_i = N-1-i+\sum_{j=i+1}^N \eps_i$, for all $1 \leq i \leq N-1$; in particular, $\beta_{N-1} = \eps_N$. For each $1 \leq i \leq N-2$,
\begin{equation}
\label{in6}
\int_0^{\pi} \cos^{\eps_i} \theta_i\ \sin^{\beta_i} \theta_i\ d\theta_i = (1+(-1)^{\eps_i})  \int_0^{\pi/2} \sin^{\eps_i} \theta_i\ \cos^{\beta_i} \theta_i\ d\theta_i = 0,
\end{equation}
unless $\eps_i$ is even. Similarly,
\begin{eqnarray}
\label{in7}
\lefteqn{\int_0^{2\pi} \cos^{\eps_{N-1}} \theta_{N-1}\ \sin^{\eps_N} \theta_{N-1}\ d\theta_{N-1}} \nonumber \\
& & = (-1)^{\eps_{N-1}+\eps_N}  \int_{-\pi}^{\pi} \cos^{\eps_{N-1}} \theta_{N-1}\ \sin^{\eps_N} \theta_{N-1}\ d\theta_{N-1} \nonumber \\
& & = (-1)^{\eps_{N-1}}(1+(-1)^{\eps_N})  \int_0^{\pi} \cos^{\eps_{N-1}} \theta_{N-1}\ \sin^{\eps_N} \theta_{N-1}\ d\theta_{N-1} \nonumber \\
& & = (-1)^{\eps_{N-1}+\eps_N}(1+(-1)^{\eps_N})  \int_{-\pi/2}^{\pi/2} \sin^{\eps_{N-1}} \theta_{N-1}\ \cos^{\eps_N} \theta_{N-1}\ d\theta_{N-1} \nonumber \\
& & = (1+(-1)^{\eps_{N-1}}+(-1)^{\eps_N}+(1)^{\eps_{N-1}+\eps_N}) \times \nonumber \\
& &\times \int_0^{\pi/2} \sin^{\eps_{N-1}} \theta_{N-1}\ \cos^{\eps_N} \theta_{N-1}\ d\theta_{N-1} = 0,
\end{eqnarray}
unless $\eps_{N-1}$ and $\eps_N$ are both even. So assume that for each $1 \leq i \leq N$, $\eps_i = 2t_i$ for some $t_i \in \zed_+$; notice that in this case the weight $w(\bm_1+\bm_2) = 2\sum_{i=1}^N t_i$, we will just call it $w$ until further notice. Then $\beta_i = N-1-i+2\sum_{j=i+1}^N t_j$, and so $(-1)^{\beta_i} = (-1)^{N-1-i}$. Putting things together, we see that $S(\bm_1,\bm_2)=0$ unless $(\bm_1,\bm_2) \in E$, in which case combining (\ref{in5}), (\ref{in6}), (\ref{in7}) and using the standard integral formulas, as for instance (5-41) and (5-42) on page 182 of \cite{fleming}, produces
\begin{eqnarray}
\label{in8}
\alpha_N S(\bm_1,\bm_2) & = & 2 \prod_{i=1}^{N-1} \left( 2 \int_0^{\pi/2} \sin^{2t_i} \theta_i\ \cos^{\beta_i} \theta_i\ d\theta_i \right) \nonumber \\
& = & 2 \left( \prod_{i=1}^{N-1} \frac{\Gamma\left( \frac{2t_i+1}{2} \right) \Gamma\left( \frac{\beta_i+1}{2} \right)}{\Gamma\left( \frac{2t_i+\beta_i}{2}+1 \right)} \right) \nonumber \\
& = & 2 \left( \prod_{i=1}^{N-1} \frac{\sqrt{\pi}  (2t_i-1)!! \Gamma\left( \frac{\beta_i+1}{2} \right)}{2^{t_i} \Gamma\left( \frac{2t_i+\beta_i}{2}+1 \right)}  \right),
\end{eqnarray}
where $\frac{2t_{i+1}+\beta_{i+1}}{2}+1 = \frac{\beta_i+1}{2}$, and so
\begin{eqnarray}
\label{in9}
\alpha_N S(\bm_1,\bm_2) & = & \frac{\Gamma\left( \frac{\beta_{N-1}+1}{2} \right)}{\Gamma\left( \frac{2t_1+\beta_1}{2}+1 \right)} 2 \pi^{\frac{N-1}{2}} \prod_{i=1}^{N-1} \frac{(2t_i-1)!!}{2^{t_i}} \nonumber \\
& = & \frac{2 \pi^{\frac{N}{2}} \prod_{i=1}^N (2t_i-1)!! }{2^{w/2} \Gamma\left( \frac{N+w}{2} \right)}  \nonumber \\
& = &  \left\{ \begin{array}{ll}
\frac{2 \pi^{\frac{N}{2}} \prod_{i=1}^N (2t_i-1)!! }{2^{w/2} \left( \frac{N - 2 + w}{2} \right) \text{\Large !}} & \mbox{if $N$ is even} \\
\frac{2 \pi^{\frac{N}{2}} \prod_{i=1}^N (2t_i-1)!! }{2^{w/2} \Gamma\left( \frac{N+w}{2} \right)} & \mbox{if $N$ is odd}
\end{array}
\right. \nonumber \\
& = &  \frac{\alpha_N \prod_{i=1}^N (\eps_i-1)!! }{\prod_{k=1}^{w/2} \left( N - 2 + 2k \right)}.
\end{eqnarray}
The result of Theorem \ref{int_lemma} now follows by combining (\ref{in1}) and (\ref{in9}).
\smallskip

\begin{rem} \label{complex} Let $K$ be a number field, and let $M(K)$ be its set of places. For each $v \in M(K)$, write $K_v$ for the completion of $K$ at $v$; in particular, if $v$ is archimedean, then $K_v = \real$ or $\cee$. Write $\Sigma^v_{N-1}$ for the unit sphere centered at $\bo$ in $K_v^N$. A standard norm on the space of homogeneous polynomials of degree $M$ in $N$ variables over $K_v$ is usually defined by
\begin{equation}
\label{norm_rem}
\|F\|_v = \left\{ \begin{array}{ll}
\left( \int_{\Sigma^v_{N-1}} |F(\bx)|_v^2 d\bx \right)^{1/2} & \mbox{if $v$ is archimedean} \\
\sup\ \left\{ |F(\bx)|_v: \bx \in \Sigma^v_{N-1} \right\} & \mbox{if $v$ is non-archimedean},
\end{array}
\right.
\end{equation}
where the measure $d\bx$ in the archimedean case is normalized so that $\int_{\Sigma^v_{N-1}} d\bx = 1$ (see, for instance, \cite{poorten:vaaler} and \cite{roy:laurent}  for details). In case $K_v=\cee$, a well-known identity (see \cite{rudin:ball}, pp. 16-17) provides
\begin{equation}
\label{in15}
\|F\|^2_v = \binom{N+M}{N}^{-1} \sum_{\bm \in \M(M,N)} \binom{M}{\bm}^{-1} |c_F(\bm)|_v^2,
\end{equation}
where $\binom{M}{\bm} = \frac{M!}{m_1!...m_N!}$. If, on the other hand, $v \nmid \infty$ then we have (see \cite{poorten:vaaler}) an identity
\begin{equation}
\label{in16}
\|F\|_v = \max_{\bm \in \M(M,N)} |c_F(\bm)|_v.
\end{equation}
Our Theorem \ref{int_lemma} can be viewed as a counterpart of these formulas when $K_v=\real$.
\end{rem}
\bigskip

\section{Proof of Theorem \ref{int_basis}}
\label{basis}

Unless stated otherwise, the notation in this section is as in the statement of Theorem \ref{int_basis}. Our argument is similar to the proof of Theorem 2.4 of \cite{me:witt}. First we recall a version of Siegel's lemma, which is essentially Theorem 2 of \cite{vaaler:siegel}.

\begin{thm} \label{vaaler} Let $V$ be as in Theorem \ref{int_basis}. Then there exists a basis $f_1,\dots,f_k$ for $V$ consisting of polynomials with integer coefficients so that
\begin{equation}
\label{siegel_bound}
\prod_{i=1}^k H(f_i) \leq H(V).
\end{equation}
\end{thm}
\smallskip

We will also need a few height-comparison lemmas. The first one is an immediate corollary of Theorem 1 of \cite{vaaler:struppeck}.

\begin{lem} \label{intersection} Let $U_1$ and $U_2$ be finite-dimensional subspaces of $\PNR$. Then
$$H(U_1 \cap U_2) \leq H(U_1) H(U_2).$$
\end{lem}

\noindent
Next let $\M = \M(d,N)$, $E = E(d,N)$, and $L = L(d,N) = |\M|$. Let us also write $\L$ for the bilinear form $\L_{d,N}$ and let $\LL = (l_{ij})_{1 \leq i,j \leq L}$ be the coefficient matrix of $\L$, so
$$\L(\bx,\bwy) = \bx^t \LL \bwy,$$
for each $\bx,\bwy \in \real^L$. Then diagonal entries of $\LL$ are equal to the corresponding coefficients of $\L$ while the off-diagonal entries are the corresponding coefficients of $\L$ multiplied by~$1/2$. Let height of $\LL$, denoted by $H(\LL)$, be the maximum of absolute values of entries of $\LL$, and let $H(\L)$ be defined in the same way as height of a polynomial is defined in (\ref{height_poly}). Then we have the following simple bound.

\begin{lem} \label{LL} Let the notation be as above, then
\begin{equation}
\label{height_L}
H(\LL) \leq H(\L) = \max_{(\bm_1,\bm_2) \in E} P \left( \frac{\bm_1+\bm_2}{2} \right) \leq  \max_{\bm \in \M} P(\bm) \leq 1.
\end{equation}
\end{lem}

\proof
Notice that
\begin{eqnarray*}
H(\LL) \leq H(\L) & = & \max_{(\bm_1,\bm_2) \in E} P \left( \frac{\bm_1+\bm_2}{2} \right) \leq \max_{\bm \in \M} P(\bm) \\
& = & \max_{\bm \in \M} \frac{\prod_{i=1}^N (2m_i-1)!! }{\prod_{k=1}^{w(\bm)} \left( N - 2 + 2k \right)} \leq \max_{\bm \in \M}  P(w(\bm),0,...,0) \\
& = & \max_{\bm \in \M} \frac{(2w(\bm)-1)!! }{\prod_{k=1}^{w(\bm)} \left( N - 2 + 2k \right)} \leq 1,
\end{eqnarray*}
which proves the lemma.
\endproof

Let us write $\|\ \|$ for the usual Euclidean norm on vectors. Let us also write $H(A) = \max_{1 \leq i,j \leq L} |a_{ij}|$ for every $L \times L$ matrix $A = (a_{ij})_{1 \leq i,j \leq L}$ with real coefficients.

\begin{lem} \label{mult} Let $f \in \PdNR$, and let $A= (a_{ij})_{1 \leq i,j \leq L}$ be an $L \times L$ real matrix. Then
$$\| \bc_f^t A \| \leq L^3 H(f) H(A).$$
In particular,
\begin{equation}
\label{lf_bound}
\| \bc_f^t \LL\| \leq L^3 H(f) H(\LL) \leq L^3 H(f),
\end{equation}
where the second inequality follows by (\ref{height_L}).
\end{lem}

\proof
Recall that our indexing set $\M = \{ \bm_1,\dots,\bm_L \}$ is lexicographically ordered, and hence 
$$\bc_f^t A = \left( \sum_{i=1}^L c_f(\bm_i) a_{i1}, ..., \sum_{i=1}^L c_f(\bm_i) a_{iL} \right).$$
Then, by Cauchy-Schwarz inequality
\begin{eqnarray*}
\| \bc_f^t A \|^2 & = &\sum_{j=1}^L\left\| \sum_{i=1}^L c_f(\bm_i) a_{ij} \right\|^2 \\
& \leq & \sum_{j=1}^L \left( \sum_{i=1}^L \| a_{ij} \|^2 \right) \left( \sum_{i=1}^L \| c_f(\bm_i) \|^2 \right) \\
& = & \| A \|^2 \| \bc_f \|^2 \leq L^6 H(A)^2 H(f)^2,
\end{eqnarray*}
where $\| A \|^2 = \sum_{i=1}^L \sum_{j=1}^L |a_{ij}|^2$, and so (\ref{lf_bound}) follows.
\endproof

The next lemma is a simple version of the well known Brill-Gordan duality principle \cite{gordan:1} (also see Theorem 1 on p. 294 of \cite{hodge:pedoe}).

\begin{lem} \label{duality} Let $f \in \PdNR$ be a polynomial with integer coefficients, and let
$$U = \{ t \in \PdNR :  \bc_f ^t A \bc_t = 0 \},$$
where $A$ is as in Lemma \ref{mult} above. Let $\gamma(f)$ be the greatest common divisor of the coordinates of the vector $\bc_f^t A$. Then
\begin{equation}
\label{dual}
H(U) = \gamma(f)^{-1} \| \bc_f^t A \|.
\end{equation}
\end{lem}
\smallskip

We will now state and prove a generalization of Theorem \ref{int_basis}.

\begin{thm} \label{gen} Let $V$ be as in Theorem \ref{int_basis}, and let $B$ be a symmetric bilinear form on $\PdNR$ with $L \times L$ coefficient matrix $\B$, meaning that
$$B(f,g) = \bc_f^t \B \bc_g,\ \forall\ f,g \in \PdNR.$$
Then there exists an orthogonal basis $g_1,\dots,g_n$ for $\left( V, B \right)$ consisting of polynomials with relatively prime integer coefficients so that
\begin{equation}
\label{gen_bound}
\prod_{i=1}^n H(g_i) \leq \left( L^3 H(\B) \right)^{\frac{n(n+1)}{2}} H(V)^n.
\end{equation}
\end{thm}

\proof
We argue by induction on $n$. First suppose that $n=1$, then pick any nonzero polynomial $g_1 \in V$ with relatively prime integer coefficients, and observe that 
$$H(g_1) = |\bc_{g_1}| \leq \| \bc_{g_1} \| = H(V),$$
where $\|\ \|$ is Euclidean norm. Next assume that $n>1$ and the theorem is true for all $1 \leq j < n$. Let $0 \neq f_1 \in V$ be a vector guaranteed by Theorem \ref{vaaler} so that
\begin{equation}
\label{s2}
H(f_1) \leq H(V)^{1/n}.
\end{equation}
First assume that $f_1$ is a non-singular point in $\left( V, B \right)$. Then
$$V_1 = \{ t \in V : B(t,f_1) = 0 \} = \{ f_1 \}^{\perp} \cap V$$
has dimension $n-1$; here 
$$\{ f_1 \}^{\perp} = \{ t \in \PdNR :  B(t,f_1) = 0 \} = \{ t \in \PdNR :  \bc_{f_1}^t \B \bc_t = 0 \}.$$
By Lemma \ref{duality}, Lemma \ref{mult}, and (\ref{s2})
\begin{equation}
\label{s2-1}
H \left( \{ f_1 \}^{\perp} \right) \leq L^3 H(\B) H(V)^{1/n}.
\end{equation}
Then by Lemma \ref{intersection} and (\ref{s2-1}) we obtain
\begin{equation}
\label{s3}
H(V_1) \leq H \left( \{ f_1 \}^{\perp} \right) H(V) \leq L^3 H(\B) H(V)^{\frac{n+1}{n}}.
\end{equation}
Since $\dim_{\real} (V_1) = n-1$, the induction hypothesis implies that there exists a basis $f_2,\dots,f_n$ for $V_1$ of polynomials with relatively prime integer coefficients such that $\B(f_i, f_j) = 0$ for all $2 \leq i \neq j \leq n$, and
\begin{eqnarray}
\label{s4} 
\prod_{i=2}^n H(f_i) & \leq & \left( L^3 H(\B) \right)^{\frac{n(n-1)}{2}} H(V_1)^{n-1} \nonumber \\
& \leq & \left( L^3 H(\B) \right)^{\frac{n^2+n-2}{2}} H(V)^{\frac{n^2-1}{n}},
\end{eqnarray}
where the last inequality follows by (\ref{s3}). Combining (\ref{s2}) and (\ref{s4}), we see that $f_1,\dots,f_n$ is a basis for $V$ satisfying (\ref{gen_bound}) so that $B(f_i,f_j) = 0$ for all $1 \leq i \neq j \leq n$.

Now assume that $f_1$ is a singular point in $\left( V, B \right)$. Since $f_1 \neq 0$,, it must be true that $c_{f_1}(\bm_j) \neq 0$ for some $\bm_j \in \M = \{ \bm_1,\dots\bm_L \}$. Let
$$T_j = \{ t \in \PdNR : c_t(\bm_j) = 0 \},$$
and define $V_1 = V \cap T_j$, then $f_1 \notin V_1$, $B(f_1,t) = 0$ for every $t \in V_1$, and 
\begin{equation}
\label{s5}
H(V_1) \leq H(V) H(T_j) = H(V),
\end{equation}
by Lemma \ref{intersection}, since $H(T_j) = 1$ by Lemma \ref{duality}. Since $\dim_{\real}(V_1) = n-1$, we can apply induction hypothesis to $V_1$, and proceed the same way as in the non-singular case above. Since the upper bound of (\ref{s5}) is smaller than that of (\ref{s3}), the result follows.
\endproof

\proof[Proof of Theorem \ref{int_basis}]
Apply Theorem \ref{gen} with $B = \L$ and use Lemma \ref{LL}.
\endproof
\bigskip

\section{Spherical designs}
\label{sphere}

In this section we apply Theorem \ref{int_lemma} to obtain a criterion for spherical designs. Let $M,N,\M(M,N)$, and $\PMNR$ be as in section \ref{intro}. A finite subset $S$ of the unit sphere $\Sigma_{N-1}$ is called a {\it spherical $M$-design} if for every $F(\bX) \in \PMNR$,
\begin{equation}
\label{des_def}
\frac{1}{\alpha_N} \int_{\Sigma_{N-1}} F(\bx)\ d\bx = \frac{1}{|S|} \sum_{\bwy \in S} F(\bwy).
\end{equation}
Spherical designs have been extensively studied, in particular in the recent years in connection with lattices, the sphere packing problem, and minimization of Epstein zeta function (see \cite{martinet} and \cite{coulangeon} for details). For instance, recent results of B. Venkov on criteria for spherical designs and their applications are summarized in Chapter 16 of \cite{martinet}. We use our Theorem \ref{int_lemma} to give another criterion for a set to be a spherical design in the general spirit of Proposition 16.1.2 and Theorem 16.1.4 of~\cite{martinet}.

\begin{thm} \label{design} Let $S \subset \Sigma_{N-1}$ be a finite set, write $S = \{ \bx_1,\dots,\bx_k \}$, where $|S|=k$. Let 
$$\M^*(M,N) = \M(M,N) \cap 2\M(M,N).$$
Then $S$ is a spherical $M$-design if and only if for every $\bm \in \M(M,N)$,
\begin{equation}
\label{crit}
 \sum_{j=1}^k \bx_j^{\bm} = \left\{ \begin{array}{ll}
k P \left( \frac{\bm}{2} \right)  & \mbox{if $\bm \in \M^*(M,N)$} \\
0 & \mbox{if $\bm \in \M(M,N) \setminus \M^*(M,N)$}.
\end{array}
\right.
\end{equation}
\end{thm}

\proof
Let us write $\bone \in \PMNR$ for the constant polynomial equal to 1, i.e. 
$$c_{\bone}(\bo) = 1,\ c_{\bone}(\bm) = 0\ \forall\ \bo \neq \bm \in \M(M,N).$$
Then for each
$$F(\bX) = \sum_{\bm \in \M(M,N)} c_F(\bm) \bX^{\bm} \in \PMNR,$$
we have
\begin{equation}
\label{d1}
\frac{1}{\alpha_N} \int_{\Sigma_{N-1}} F(\bx)\ d\bx = \left< F,\bone \right> = \L_{M,N}(\bc_F,\bc_{\bone})  = \sum_{\bm \in \M^*(M,N)} P \left( \frac{\bm}{2} \right) c_F(\bm).
\end{equation}
On the other hand, (\ref{des_def}) implies that $S$ is a spherical $M$-design if and only if
\begin{eqnarray}
\label{d2}
\frac{1}{\alpha_N} \int_{\Sigma_{N-1}} F(\bx)\ d\bx & = & \frac{1}{k} \sum_{j=1}^k \sum_{\bm \in \M(M,N)} c_F(\bm) \bx_j^{\bm} \nonumber \\
& = & \sum_{\bm \in \M(M,N)} \left( \frac{1}{k} \sum_{j=1}^k \bx_j^{\bm} \right) c_F(\bm).
\end{eqnarray}
Hence we must have
\begin{equation}
\label{des}
\sum_{\bm \in \M^*(M,N)} P \left( \frac{\bm}{2} \right) c_F(\bm) = \sum_{\bm \in \M(M,N)} \left( \frac{1}{k} \sum_{j=1}^k \bx_j^{\bm} \right) c_F(\bm),
\end{equation}
for all $F(\bX) \in \PMNR$, which means that the linear forms in the variables $\bc_F$ on the right and left hand sides of (\ref{des}) must be equal identically, i.e. their respective coefficients must be equal. Then (\ref{crit}) follows.
\endproof
\bigskip

%\nocite{*}
\bibliographystyle{plain}  % Here the bibliography 
\bibliography{integral}        % is inserted.

\end{document}